\documentclass[11pt]{amsart}
\usepackage[utf8]{inputenc}
\usepackage[T1]{fontenc}
\usepackage{amsmath,amssymb,amsthm}
\usepackage{booktabs}
\usepackage[margin=1.15in]{geometry}
\usepackage[hidelinks]{hyperref}

\newtheorem{theorem}{Theorem}[section]
\newtheorem{lemma}[theorem]{Lemma}
\newtheorem{proposition}[theorem]{Proposition}
\newtheorem{corollary}[theorem]{Corollary}
\theoremstyle{definition}
\newtheorem{definition}[theorem]{Definition}
\theoremstyle{remark}
\newtheorem{remark}[theorem]{Remark}

\newcommand{\dn}{{\downarrow}}
\newcommand{\rk}{\operatorname{rk}}
\newcommand{\Aut}{\operatorname{Aut}}
\newcommand{\Sym}{\operatorname{Sym}}
\newcommand{\N}{\mathbb{N}}
\newcommand{\Ndist}{160}
\newcommand{\Nmod}{38}
\newcommand{\Nmodprev}{37}
\newcommand{\Ndistprev}{159}

\newcommand{\ratiod}{1.84337061669}
\newcommand{\ratiov}{1.72654}
\newcommand{\ratiom}{2.4194}
\newcommand{\ratiomv}{2.4072}
\newcommand{\tabledist}{\begin{tabular}{rrr}\toprule $n$ & $d(n)$ & $d_{\mathrm{vi}}(n)$\\\midrule
61 & {\footnotesize $1\,398\,297\,238\,937\,456$} & {\footnotesize $3\,563\,549\,771\,088$}\\
62 & {\footnotesize $2\,577\,577\,210\,993\,326$} & {\footnotesize $6\,180\,294\,770\,133$}\\
63 & {\footnotesize $4\,751\,419\,899\,674\,101$} & {\footnotesize $10\,627\,703\,640\,278$}\\
64 & {\footnotesize $8\,758\,621\,005\,139\,604$} & {\footnotesize $18\,415\,015\,847\,260$}\\
65 & {\footnotesize $16\,145\,360\,828\,911\,176$} & {\footnotesize $31\,692\,379\,806\,500$}\\
66 & {\footnotesize $29\,761\,867\,322\,814\,983$} & {\footnotesize $54\,874\,841\,310\,053$}\\
67 & {\footnotesize $54\,862\,096\,251\,848\,441$} & {\footnotesize $94\,501\,198\,486\,814$}\\
68 & {\footnotesize $101\,131\,136\,719\,043\,728$} & {\footnotesize $163\,532\,825\,001\,695$}\\
69 & {\footnotesize $186\,422\,036\,415\,293\,546$} & {\footnotesize $281\,769\,158\,221\,666$}\\
70 & {\footnotesize $343\,644\,809\,421\,288\,784$} & {\footnotesize $487\,372\,342\,417\,275$}\\
71 & {\footnotesize $633\,464\,442\,109\,685\,434$} & {\footnotesize $840\,095\,279\,546\,933$}\\
72 & {\footnotesize $1\,167\,709\,511\,857\,902\,566$} & {\footnotesize $1\,452\,567\,011\,369\,089$}\\
73 & {\footnotesize $2\,152\,520\,697\,500\,238\,180$} & {\footnotesize $2\,504\,648\,622\,120\,683$}\\
74 & {\footnotesize $3\,967\,892\,860\,459\,105\,208$} & {\footnotesize $4\,329\,391\,888\,717\,253$}\\
75 & {\footnotesize $7\,314\,295\,461\,557\,816\,868$} & {\footnotesize $7\,467\,094\,119\,515\,045$}\\
76 & {\footnotesize $13\,482\,956\,030\,171\,180\,202$} & {\footnotesize $12\,904\,165\,237\,004\,308$}\\
77 & {\footnotesize $24\,854\,081\,123\,501\,649\,639$} & {\footnotesize $22\,261\,052\,808\,306\,244$}\\
78 & {\footnotesize $45\,815\,279\,724\,188\,878\,657$} & {\footnotesize $38\,462\,962\,752\,711\,714$}\\
79 & {\footnotesize $84\,454\,531\,446\,767\,120\,345$} & {\footnotesize $66\,363\,791\,347\,778\,827$}\\
80 & {\footnotesize $155\,680\,994\,246\,494\,614\,800$} & {\footnotesize $114\,647\,177\,779\,433\,135$}\\
90 & {\footnotesize $70\,528\,082\,268\,645\,205\,822\,241$} & {\footnotesize $26\,979\,773\,850\,610\,228\,676$}\\
100 & {\footnotesize $31\,951\,309\,377\,980\,462\,944\,738\,960$} & {\footnotesize $6\,349\,929\,081\,499\,442\,762\,100$}\\
110 & {\footnotesize $14\,474\,889\,691\,239\,516\,050\,456\,290\,105$} & {\footnotesize $1\,494\,573\,995\,405\,456\,856\,950\,573$}\\
120 & {\footnotesize $6\,557\,553\,883\,147\,316\,182\,319\,056\,599\,661$} & {\footnotesize $351\,780\,516\,810\,470\,934\,005\,984\,975$}\\
130 & {\footnotesize $2\,970\,766\,195\,151\,439\,931\,838\,036\,667\,770\,108$} & {\footnotesize $82\,799\,559\,322\,022\,867\,662\,771\,980\,605$}\\
140 & {\footnotesize $1\,345\,845\,103\,971\,103\,748\,656\,962\,514\,254\,560\,619$} & {\footnotesize $19\,488\,790\,063\,093\,416\,469\,215\,478\,130\,947$}\\
150 & {\footnotesize $609\,707\,706\,680\,280\,574\,986\,503\,202\,459\,978\,289\,571$} & {\footnotesize $4\,587\,139\,276\,311\,270\,850\,485\,527\,015\,829\,145$}\\
160 & {\footnotesize $276\,215\,655\,494\,051\,950\,398\,467\,067\,725\,855\,360\,399\,901$} & {\footnotesize $1\,079\,689\,902\,088\,453\,963\,903\,370\,538\,000\,719\,118$}\\
\bottomrule\end{tabular}}
\newcommand{\tablemod}{\begin{tabular}{rrr}\toprule $n$ & $m(n)$ & $m_{\mathrm{vi}}(n)$\\\midrule
33 & $48\,496\,079\,939$ & $6\,579\,923\,491$\\
34 & $116\,905\,715\,114$ & $15\,777\,658\,535$\\
35 & $282\,098\,869\,730$ & $37\,871\,501\,929$\\
36 & $681\,357\,605\,302$ & $90\,998\,884\,153$\\
37 & $1\,647\,135\,247\,659$ & $218\,856\,768\,070$\\
38 & $3\,985\,106\,742\,170$ & $526\,836\,817\,969$\\
\bottomrule\end{tabular}}
\newcommand{\resultsother}{\begin{table}[ht]
\begin{tabular}{rr}\toprule $n$ & semimodular\\\midrule
23 & $312\,324\,027$\\
24 & $1\,086\,545\,705$\\
25 & $3\,838\,581\,926$\\
\midrule
26 & $13\,759\,007\,833$\\
27 & $49\,997\,352\,868$\\
\bottomrule\end{tabular}\qquad
\caption{Number of semimodular lattices with $n$ elements (A229202).  The values above the rule agree with~\cite{Koh19}; those below are new.}\label{tab:other}
\end{table}}
\newcommand{\boundedtable}{\begin{table}[ht]
\begin{tabular}{rrrll}\toprule $W$ & states & edges & growth rate, v.i. & growth rate, all\\\midrule
2 & $3$ & $7$ & $1.4142135624$ & $1.7106440950$\\
3 & $13$ & $52$ & $1.6616501964$ & $1.8211213614$\\
4 & $25$ & $143$ & $1.7168846567$ & $1.8411147093$\\
5 & $50$ & $351$ & $1.7235737135$ & $1.8425913201$\\
6 & $266$ & $2\,327$ & $1.7262480148$ & $1.8433325496$\\
7 & $735$ & $7\,031$ & $1.7264936580$ & $1.8433671106$\\
8 & $2\,203$ & $23\,668$ & $1.7265344528$ & $1.8433705132$\\
9 & $6\,616$ & $70\,646$ & $1.7265353517$ & $1.8433705490$\\
10 & $22\,620$ & $263\,509$ & $1.7265362224$ & $1.8433706134$\\
11 & $81\,606$ & $971\,174$ & $1.7265363005$ & $1.8433706166$\\
12 & $274\,115$ & $3\,309\,503$ & $1.7265363021$ & $1.8433706167$\\
\bottomrule\end{tabular}
\caption{The automata $\mathcal{A}_W$ and the exponential growth rates of the distributive lattices all of whose layers have at most $W$ elements (floating point values; the certified bounds are in Theorem~\ref{thm:growth}).}\label{tab:bounded}
\end{table}}
\newcommand{\boundedtablemod}{\begin{table}[ht]
\begin{tabular}{rrrll}\toprule $W$ & states & edges & growth rate, v.i. & growth rate, all\\\midrule
2 & $3$ & $7$ & $1.4142135624$ & $1.7106440950$\\
3 & $13$ & $80$ & $1.9139253861$ & $2.0654379980$\\
4 & $41$ & $533$ & $2.1405885481$ & $2.2291160113$\\
5 & $139$ & $4\,036$ & $2.2583196481$ & $2.3158654358$\\
6 & $785$ & $45\,457$ & $2.3262518513$ & $2.3664739489$\\
7 & $5\,359$ & $839\,353$ & $2.3680771954$ & $2.3978761910$\\
8 & $31\,435$ & $10\,654\,239$ & $2.3951359778$ & $2.4182279039$\\
\bottomrule\end{tabular}
\caption{The automata $\mathcal{A}^{\mathrm{mod}}_W$ and the exponential growth rates of the modular lattices all of whose layers have at most $W$ elements (floating point values; certified bounds in Theorem~\ref{thm:modgrowth}).}\label{tab:boundedmod}
\end{table}}

\title[Counting distributive and modular lattices]{Counting distributive and modular lattices\\ with a frontier transfer matrix}
\author{Jobst Heitzig}
\address{Potsdam Institute for Climate Impact Research, Telegrafenberg A31, 14473 Potsdam, Germany}

\begin{document}

\begin{abstract}
We describe a transfer-matrix method for counting unlabelled graded lattices by
their rank layers, and apply it to distributive and to modular lattices.  Earlier
counts generate the lattices one at a time.  Here all lower parts of lattices
that look the same from above are merged into one \emph{state}.  For a
distributive lattice $J(P)$ the state after rank $r$ is the poset of those elements of
$P$ that lie in some but not all down-sets of size $r$, together with the
symmetry group induced on it; for a modular lattice it is the poset of pairwise
meets of the elements of rank $r$.  We prove that these states determine all
possible continuations up to isomorphism, so that the work grows with the number
of states and not with the number of lattices.  For modular lattices the number
of states is reduced further by a bound that rests on self-duality of the class.
We obtain the number of distributive lattices up to $\Ndist$ elements; the previous record
was $60$.  For modular lattices we reproduce the known numbers up to $\Nmod$ elements in
$142$ minutes on one processor core.
The frontier for modular lattices also serves for semimodular lattices, for which we
obtain two further values ($26$ and $27$ elements).
For distributive lattices with layers of bounded size the transfer matrix is a finite
automaton; from it we derive the certified bounds $1.7265363$ and $1.84337061$ for the
exponential growth rates of the vertically indecomposable and of all distributive lattices
(previously $1.7250$ and $1.84$), and $2.3951$ for vertically indecomposable modular
lattices (previously $2.3122$).
The $60$ previously known values for distributive lattices take $0.1$ seconds.  The
new values are supported by the known ones, by duality checks, by brute force
over all posets on at most $11$ points, and by agreement of the two programs on
distributive lattices.
\end{abstract}

\maketitle

\section{Introduction}

Let $d(n)$ and $m(n)$ be the numbers of isomorphism classes of distributive and of
modular lattices with $n$ elements (sequences A006982 and A006981 in the
OEIS~\cite{oeis}).  A lattice is \emph{vertically indecomposable} if no element other
than the top and bottom is comparable to all elements; $d_{\mathrm{vi}}(n)$ and
$m_{\mathrm{vi}}(n)$ count such lattices with $n \ge 2$ elements (A072361, A342132).
Every lattice is in a unique way a vertical sum of vertically indecomposable ones, so
\begin{equation}\label{eq:vsum}
  d(n) = \sum_{k=2}^{n} d_{\mathrm{vi}}(k)\, d(n-k+1) \quad (n \ge 2), \qquad d(1)=1,
\end{equation}
and likewise for $m$.

\subsection*{Previous work}
Ern\'e, Heitzig and Reinhold~\cite{EHR} determined $d(n)$ for $n \le 49$ by generating
posets.  Kohonen~\cite{Koh21} reached $n \le 60$ by generating only those lattices
that are not vertical $2$-sums of smaller ones; the computation for $n=60$ alone took
$17.2$ core-days, with a factor of about $1.66$ per additional element.
Modular lattices were counted for $n \le 24$ by Jipsen and Lawless~\cite{JL} (the OEIS records a later correction of the value for $n=24$), for
$n \le 30$ by Kohonen~\cite{Koh19}, and for $n \le 35$ in~\cite{Koh21}, where $n=35$
took $63.1$ core-days, with a factor of about $2.66$ per element.  In~\cite{Koh23} the
count was extended to $n \le 40$ by listing only the \emph{racks}, i.e.\ the lattices from
which certain doubly irreducible elements have been removed, and counting the remaining
lattices by P\'olya's method; listing the racks took $70.4$ core-days, of which $8.1$ for
$n = 38$ and $38.5$ for $n = 40$.  All these
methods construct the lattices (or a large part of them) individually and remove
isomorphic copies with canonical labellings~\cite{nauty}; see
also~\cite{HR,GT} for general lattices and~\cite{BM} for posets.

\subsection*{This paper}
We do not generate lattices.  A finite distributive or modular lattice is graded, and we
count by a transfer matrix over its rank layers, in the spirit of the transfer-matrix
method for walks in a graph~\cite[Section~4.7]{Stanley}.  The difficulty is to find
a \emph{state} that (a) determines everything that can happen above a layer, (b) is
small for most lattices, and (c) allows isomorphism classes, rather than labelled
objects, to be counted.  Sections~\ref{sec:dist} and~\ref{sec:mod} give such states.
Point (c) is handled by making the image of the automorphism group of the lower part
a component of the state (Lemmas~\ref{lem:sym} and~\ref{lem:modsym}).

The lattices with few elements per layer are the vast majority, and for them the states are
tiny.  For instance all distributive lattices with at most $60$ elements, about
$1.7\cdot 10^{15}$ of them, are accounted for by $4496$ states.

\subsection*{Results}
\begin{itemize}
\item $d(n)$ and $d_{\mathrm{vi}}(n)$ for all $n \le \Ndist$ (Table~\ref{tab:dist});
  the first new value is\\ $d(61) = 1\,398\,297\,238\,937\,456$.
\item $m(n)$ and $m_{\mathrm{vi}}(n)$ for all $n \le \Nmod$ (Table~\ref{tab:mod}).  These are
  not new: they agree with~\cite{Koh23}, where $n \le 40$ is reached.  What is new is the
  method and its cost.
\item Certified lower bounds for the exponential growth rates of distributive and of modular
  lattices (Sections~\ref{sec:bounded} and~\ref{sec:mod}).
\item Two sequences that are functions of $d$ and $d_{\mathrm{vi}}$ (Section~\ref{sec:further}).
\item The numbers of semimodular lattices with $26$ and $27$ elements
  (Section~\ref{sec:other}); the previous record was $25$~\cite{Koh19}.
\end{itemize}
The complete tables, the programs and the logs of all runs are available as
ancillary files and at~\cite{code}.  The values are the output of
programs, and Section~\ref{sec:verify} describes the evidence for them.

\section{Distributive lattices}\label{sec:dist}

By Birkhoff's theorem~\cite[Theorem~3.4.1]{Stanley} every finite distributive lattice
is isomorphic to the lattice $J(P)$ of down-sets of a finite poset $P$, unique up to
isomorphism.  So $d(n)$ is the number of unlabelled posets with exactly $n$ down-sets.
The rank of a down-set in $J(P)$ is its cardinality.

\subsection{Layers, births and the frontier}
Let $P$ be a finite poset.  For $x \in P$ put $h(x) = |\dn x|$, where
$\dn x = \{y : y \le x\}$.  For $0 \le r \le |P|$ let
\[
  X_r = \{x \in P : h(x) \le r\}, \qquad S_r = \{I \in J(P) : |I| = r\}.
\]
$X_r$ is a down-set of $P$, since $y < x$ implies $h(y) < h(x)$.  Every $I \in S_r$ is
contained in $X_r$, because $\dn x \subseteq I$ for $x \in I$.  Conversely every
$x \in X_r$ lies in some $I \in S_r$: extend $\dn x$ by minimal elements of the
complement until it has $r$ elements.  We say that $x$ is \emph{born at rank} $h(x)-1$;
its strict down-set $\dn x \setminus \{x\}$ belongs to $S_{h(x)-1}$.

\begin{definition}
An element of $X_r$ is \emph{dead} (at rank $r$) if it lies in every member of
$S_r$, and \emph{live} otherwise.  Let $D_r$ be the set of dead elements,
$\Lambda_r$ the subposet of $P$ induced on the live elements, and
$\rho_r = r - |D_r|$.  The pair $(\Lambda_r,\rho_r)$ is the \emph{frontier} of $P$ at
rank $r$.
\end{definition}

\begin{lemma}\label{lem:frontier}
\textup{(a)} $D_r$ is a down-set of $P$, and $D_r \subseteq D_{r+1}$.

\textup{(b)} $I \mapsto I \setminus D_r$ is a bijection from $S_r$ onto the set of
down-sets of $\Lambda_r$ of size $\rho_r$.

\textup{(c)} If $h(x) = r+1$, then $\dn x \setminus\{x\} = D_r \cup K$ for a down-set $K$ of $\Lambda_r$
of size $\rho_r$, and $x$ lies above every dead element.
\end{lemma}

\begin{proof}
(a) If $x \in D_r$ and $y < x$, then $y$ lies in every down-set containing $x$.  Every
down-set of size $r+1$ contains one of size $r$ (remove a maximal element), hence
contains $D_r$.

(b) Let $I \in S_r$.  Then $D_r \subseteq I$, and $I \setminus D_r$ is a down-set of
$\Lambda_r$ with $\rho_r$ elements.  Conversely let $K$ be a down-set of $\Lambda_r$
with $|K| = \rho_r$.  An element below a member of $D_r \cup K$ lies in $X_r$; it is
dead, or it is live and then belongs to $K$.  So $D_r \cup K \in S_r$.

(c) follows from (b) because $\dn x \setminus \{x\} \in S_r$.
\end{proof}

Thus the elements born at rank $r$ are described by a \emph{birth vector}
\[
  \mu \colon \{\text{down-sets of } \Lambda_r \text{ of size } \rho_r\} \to \N,
  \qquad \mu(K) = \#\{x : \dn x \setminus\{x\} = D_r \cup K\}.
\]

\begin{lemma}\label{lem:next}
The frontier $(\Lambda_{r+1},\rho_{r+1})$ and the number $|S_{r+1}|$ are determined by
$(\Lambda_r,\rho_r)$ and $\mu$.
\end{lemma}

\begin{proof}
Let $\Lambda'$ be $\Lambda_r$ together with the elements born at rank $r$; these are
new maximal elements, pairwise incomparable, with down-sets given by $\mu$.  By
Lemma~\ref{lem:frontier}(a) every member of $S_{r+1}$ contains $D_r$, and as in the
proof of Lemma~\ref{lem:frontier}(b) the map $I \mapsto I \setminus D_r$ is a
bijection from $S_{r+1}$ to the down-sets of $\Lambda'$ of size $\rho_r + 1$.  Their
intersection $E$ gives $D_{r+1} = D_r \cup E$, so $\Lambda_{r+1} = \Lambda' \setminus E$
and $\rho_{r+1} = \rho_r + 1 - |E|$.
\end{proof}

So the part of $J(P)$ above rank $r$ depends on $X_r$ only through the frontier.  The
frontier is empty exactly when $|S_r| = 1$, that is, when $J(P)$ has a single element of
rank $r$.

\subsection{Isomorphism classes}
A \emph{lower part of height $r$} is a finite poset $X$ with $h(x) \le r$ for all
$x$ and $|X| \ge r$.  For any poset $P$ and $r \le |P|$ the set $X_r$ is one, and $P$ is
a lower part of height $|P|$.  The \emph{weight} of $X$ is
$\omega(X) = \sum_{i=0}^{r} |S_i(X)|$; note that $\omega(P) = |J(P)|$.  An
\emph{extension} of $X$ is a lower part $X'$ of height $r+1$ with $X_r(X') = X$; it is
determined up to isomorphism over $X$ by its birth vector.

Since $h$ is invariant under isomorphisms, an isomorphism between posets maps $X_r$ to
$X_r$, $S_r$ to $S_r$, and dead elements to dead elements.  In particular $\Aut(X)$ acts on
$\Lambda_r(X)$; let $\Gamma(X) \le \Aut(\Lambda_r(X))$ be the image of this action.  The group
$\Gamma(X)$ permutes the down-sets of $\Lambda_r$ of size $\rho_r$ and hence the birth
vectors.

\begin{lemma}\label{lem:sym}
Let $X$ be a lower part of height $r$, and let $X'$, $X''$ be the extensions with birth
vectors $\mu'$, $\mu''$.

\textup{(a)} $X' \cong X''$ if and only if $\mu'' = \mu' \circ \gamma$ for some
$\gamma \in \Gamma(X)$.

\textup{(b)} $\Gamma(X')$ is determined by $(\Lambda_r(X),\rho_r(X),\Gamma(X))$ and
$\mu'$.
\end{lemma}

\begin{proof}
(a) An isomorphism $X' \to X''$ restricts to an automorphism $\alpha$ of
$X = X_r(X') = X_r(X'')$, and it maps the elements born over $D_r \cup K$ to the
elements born over $D_r \cup \alpha(K)$; so $\mu'(K) = \mu''(\alpha K)$.  Conversely, if
$\mu'' \circ \alpha|_{\Lambda_r} = \mu'$, then $\alpha$ extends to an isomorphism by any
bijection of the new elements that is compatible with $\alpha$.

(b) By the same argument $\Aut(X')$ consists of the pairs $(\alpha,\beta)$ where
$\alpha \in \Aut(X)$ fixes $\mu'$ and $\beta$ is a compatible bijection of the new
elements.  Its action on $\Lambda_{r+1}(X') \subseteq \Lambda_r(X) \cup \{\text{new
elements}\}$ only involves the image of $\alpha$ in $\Gamma(X)$ and $\beta$.
\end{proof}

Call two triples $(\Lambda,\rho,\Gamma)$, $(\Lambda',\rho,\Gamma')$ with
$\Gamma \le \Aut(\Lambda)$ equivalent if a poset isomorphism $\Lambda \to \Lambda'$
carries $\Gamma$ to $\Gamma'$; a \emph{state} is an equivalence class.  The state of a
lower part $X$ of height $r$ is the class of
$(\Lambda_r(X),\rho_r(X),\Gamma(X))$.  For a state $\sigma$ let $w(\sigma)$ be the
number of down-sets of $\Lambda$ of size $\rho$, and for a $\Gamma$-orbit $[\mu]$ of
birth vectors let $\sigma \cdot [\mu]$ be the state of the extension, which is well
defined by Lemmas~\ref{lem:next} and~\ref{lem:sym}.  Write $\varepsilon$ for the state
with empty $\Lambda$.

\begin{theorem}\label{thm:dist}
Let $c_r(\sigma,a)$ be the number of isomorphism classes of lower parts of height $r$
with state $\sigma$ and weight $a$.  Then $c_0(\varepsilon,1)=1$ and
\[
  c_{r+1}(\tau,a) \;=\; \sum_{\sigma}\ \sum_{[\mu]\,:\,\sigma\cdot[\mu]=\tau}
     c_r\bigl(\sigma,\,a - w(\tau)\bigr),
\]
where $[\mu]$ runs over the $\Gamma$-orbits of birth vectors of $\sigma$ for which the
extension has a down-set of size $r+1$.  Moreover
\[
  d(n) = \sum_{r \ge 0} c_r(\varepsilon,n) ,
\]
and $d_{\mathrm{vi}}(n)$ is obtained in the same way if all lower parts that have state
$\varepsilon$ at some height strictly between $0$ and their own are left out.
\end{theorem}

\begin{proof}
The map $[X'] \mapsto [X_r(X')]$ from isomorphism classes of lower parts of height
$r+1$ to those of height $r$ is well defined, and by Lemma~\ref{lem:sym}(a) its fibre
over $[X]$ is in bijection with the $\Gamma(X)$-orbits of admissible birth vectors.  The
weight increases by $|S_{r+1}(X')| = w(\tau)$.  A lower part of height $r$ has state
$\varepsilon$ if and only if it has exactly $r$ elements, i.e.\ it is a poset $P$ with
$|P| = r$ regarded as its own lower part; its weight is $|J(P)|$.  Finally $J(P)$ is
vertically indecomposable if and only if $|S_i| \ge 2$ for $0 < i < |P|$.
\end{proof}

\begin{proposition}\label{prop:bound}
Let $X$ be a lower part with frontier $(\Lambda,\rho)$, and let $u$ be the number of
down-sets of $\Lambda$ of size greater than $\rho$.  Every poset $P$ with $X_r(P) = X$
has $|J(P)| \ge \omega(X) + u$, with equality for $P = X$.
\end{proposition}

\begin{proof}
For each such down-set $K$ the set $D_r \cup K$ is a down-set of $X$, hence of $P$, of
size greater than $r$.  If $P = X$, every down-set of size greater than $r$ contains a
member of $S_r$, hence $D_r$, and is of this form.
\end{proof}

So a state reached with weight $a$ needs to be kept only if $a + u \le N$, where $N$
is the largest $n$ wanted, and this test is exact: every state that is kept occurs in
some lattice with at most $N$ elements.

\subsection{Implementation}
States are stored in a hash table under a canonical form: $\Lambda$ is relabelled
canonically (by iterated hashing of invariants, with \texttt{nauty}~\cite{nauty} as
fallback when the resulting partition is not discrete), and $\Gamma$ is stored as the
sorted list of its elements in these labels.  We do not normalise $\Gamma$ under
conjugation by $\Aut(\Lambda)$; this may store one state several times but does not
affect correctness.  The counters $c(\sigma,a)$ are processed in increasing order of
$a$ (the height $r$ need not be recorded).  For each state the birth vectors are
enumerated by depth-first search, pruned with Proposition~\ref{prop:bound}, and one
representative per $\Gamma$-orbit is kept (the lexicographically least).  Twins,
i.e.\ several elements born over the same down-set, contribute a symmetric group to
$\Gamma$; in a distributive lattice their number is small, because $t$ twins force
$2^t$ down-sets.  Counters have $160$ bits and overflow is trapped.

\subsection{Results}
Table~\ref{tab:cost} shows the cost on one core of an Intel Xeon at $2.3$\,GHz.  The
number of states grows by a factor of about $2.1$ when $N$ increases by $10$, while
$d(N)$ grows by a factor of about $450$.

\begin{table}[ht]
\begin{tabular}{rrrr}\toprule
$N$ & states & time & memory\\\midrule
$60$ & $4.5\cdot 10^{3}$ & $0.1$ s & \\
$100$ & $1.9\cdot 10^{5}$ & $9$ s & \\
$120$ & $9.7\cdot 10^{5}$ & $41$ s & $0.4$ GB\\
$150$ & $8.8\cdot 10^{6}$ & $7.5$ min & $3.9$ GB\\
$160$ & $1.75\cdot 10^{7}$ & $16$ min & $7.9$ GB\\\bottomrule
\end{tabular}
\caption{Cost of computing $d_{\mathrm{vi}}(n)$ for all $n \le N$.}\label{tab:cost}
\end{table}

\begin{table}[ht]
\centering\setlength{\tabcolsep}{4pt}
\tabledist
\caption{Values of $d(n)$ (A006982) and $d_{\mathrm{vi}}(n)$ (A072361) for $n > 60$.  The
computed values for $n \le 60$ agree with~\cite{Koh21}.  All values for $n \le \Ndist$ are in the ancillary files.}\label{tab:dist}
\end{table}

\subsection{Growth}
Let $F(x) = \sum_{n \ge 1} d(n) x^{n-1}$ and $V(x) = \sum_{k \ge 2} d_{\mathrm{vi}}(k) x^{k-1}$.
Equation~\eqref{eq:vsum} reads $F = 1/(1-V)$.

\begin{proposition}\label{prop:growth}
\textup{(a)} The limit $\beta = \lim_{n\to\infty} d(n)^{1/n}$ exists and is finite.

\textup{(b)} If $\limsup_{n} d_{\mathrm{vi}}(n)^{1/n} < \beta$, then
$d(n) \sim C \beta^n$ with $C = 1/V'(1/\beta)$.
\end{proposition}

\begin{proof}
(a) This is a special case of~\cite[Corollary~11]{EHR}, which applies to every sequence
defined by a recursion of the form~\eqref{eq:vsum}.  A direct argument: the vertical sum of
lattices with $i$ and $j$ elements has $i+j-1$ elements, and
the two summands can be recovered from it and $i$.  So $d(i+j-1) \ge d(i)\,d(j)$, the
sequence $d(n+1)$ is supermultiplicative, and the limit exists by Fekete's lemma.  It
is finite by the exponential upper bound of~\cite{EHR}.

(b) Put $\varrho = 1/\beta$.  By assumption $V$ is analytic in a disc of radius larger
than $\varrho$, while $F$ has radius of convergence $\varrho$ and, having nonnegative
coefficients, a singularity at $\varrho$.  Hence $V(\varrho) = 1$.  Since
$V'(\varrho) > 0$ this is a simple pole of $F$, and it is the only singularity on
$|x| = \varrho$: $V(x) = 1$ with $|x| = \varrho$ forces all terms of $V(x)$ to be
positive, and the term $d_{\mathrm{vi}}(2)\,x$ is present.  The claim follows by the
standard transfer from a dominant simple pole~\cite[Section~V.2]{FS}.
\end{proof}

Part (b) is the standard behaviour of a supercritical sequence construction and is not
new either; we state it to make precise what the data suggest.  By~\cite[Corollary~11]{EHR}
and the canonical $2$-decomposition of~\cite{EHR} the limit of $d_{\mathrm{vi}}(n)^{1/n}$
exists as well.  The hypothesis of (b) is open.  The data support it strongly:
$d(\Ndist)/d(\Ndistprev) = \ratiod$, which agrees with the ratio at $n=150$ to twelve
digits, while $d_{\mathrm{vi}}(\Ndist)/d_{\mathrm{vi}}(\Ndistprev) = \ratiov$.
Numerically $\beta \approx 1.84337061669$ and $C \approx 0.087745045$.
Kohonen~\cite{Koh21} proved $d_{\mathrm{vi}}(n) = \Omega(1.7250^n)$ and remarked that
this bound appears to be close to the truth; the ratio $\ratiov$ supports that.  A
proof that vertically indecomposable distributive lattices have a strictly smaller
exponential growth rate than all distributive lattices would turn the numerical
statement into a theorem.  Rigorous lower bounds for both growth rates are given in
Section~\ref{sec:bounded}.

\section{Bounded layers and growth rates}\label{sec:bounded}

Let $\lambda_{\mathrm{vi}} = \limsup_n d_{\mathrm{vi}}(n)^{1/n}$ (by~\cite[Corollaries~7 and~11]{EHR}
this is in fact a limit), and let
$\beta$ be as in Proposition~\ref{prop:growth}.  The best published bounds are
$\lambda_{\mathrm{vi}} \ge 1.7250$~\cite{Koh21} and $\beta \ge 1.84$~\cite{EHR}.  The transfer
matrix of Section~\ref{sec:dist} gives much sharper lower bounds, because for lattices with
layers of bounded size it is a \emph{finite} automaton.

For $W \ge 2$ let $\mathcal{D}_W$ be the class of finite distributive lattices all of whose
rank layers have at most $W$ elements.  A lower part belongs to a lattice of $\mathcal{D}_W$
only if every layer that is already forced has at most $W$ elements, i.e.\ if its
frontier $(\Lambda,\rho)$ has at most $W$ down-sets of each size $\ge \rho$
(Proposition~\ref{prop:bound}).  Let $\mathcal{A}_W$ be the directed multigraph whose
vertices are the states with this property that can be reached from $\varepsilon$, with one
edge for every $\Gamma$-orbit of birth vectors, and with the weight $w(\sigma) \le W$ on the
vertex $\sigma$.  By Theorem~\ref{thm:dist} the lattices in $\mathcal{D}_W$ with $n$ elements
correspond to the closed walks at $\varepsilon$ whose entered vertices have total weight $n-1$,
and the vertically indecomposable ones to the walks that do not pass through $\varepsilon$.

\begin{proposition}\label{prop:finite}
For $2 \le W \le 12$ the multigraph $\mathcal{A}_W$ is finite, with the numbers of vertices
given in Table~\ref{tab:bounded}.  Consequently the generating function of $\mathcal{D}_W$ is
rational for these $W$.
\end{proposition}

This is the outcome of a computation (the exploration of the reachable states terminates); we
have no proof for general $W$.

\begin{lemma}\label{lem:cw}
Let $C$ be a strongly connected set of vertices of a finite directed multigraph with positive
integer vertex weights $w$ and $e(s,t)$ edges from $s$ to $t$.  Suppose that $x > 0$ and
$u \colon C \to \mathbb{R}_{>0}$ satisfy
\begin{equation}\label{eq:cw}
  \sum_{t \in C} e(s,t)\, x^{w(t)}\, u_t \;\ge\; u_s \qquad\text{for all } s \in C .
\end{equation}
Then for every $s \in C$ the number $c_s(n)$ of closed walks at $s$ inside $C$ whose entered
vertices have total weight $n$ satisfies $\limsup_n c_s(n)^{1/n} \ge 1/x$.
\end{lemma}

\begin{proof}
Let $M = (e(s,t)\,x^{w(t)})_{s,t \in C}$.  Then $\sum_n c_s(n) x^n = \sum_k (M^k)_{ss}$.  By
\eqref{eq:cw}, $M^k u \ge u$ for all $k$.  If $\sum_k (M^k)_{ss}$ converged for one $s$, then,
$M$ being irreducible, $\sum_k (M^k)_{tt'}$ would converge for all $t,t'$ (each entry of
$M^k$ is bounded by a constant times a diagonal entry of $M^{k+j}$ for a fixed $j$), and
$M^k u \to 0$, a contradiction.  So the power series $\sum_n c_s(n)x^n$ diverges at $x$.
\end{proof}

\begin{theorem}\label{thm:growth}
$\lambda_{\mathrm{vi}} \ge 1.7265363$ and $\beta \ge 1.84337061$.
\end{theorem}

\begin{proof}
Take $W = 12$.  For the first bound let $C$ be the strong component of
$\mathcal{A}_{12} \setminus \{\varepsilon\}$ with $274\,096$ vertices; it can be reached from
$\varepsilon$ and $\varepsilon$ can be reached from it.  Different walks give non-isomorphic
lattices, so $d_{\mathrm{vi}}(n + j) \ge c_s(n)$ for a fixed $j$ and $s \in C$.  For the second
bound let $C$ be all of $\mathcal{A}_{12}$, which is strongly connected, and $s = \varepsilon$;
then $d(n+1) \ge c_\varepsilon(n)$.  In both cases a positive integer vector $u$ satisfying
\eqref{eq:cw} with $x = 1/1.7265363$, respectively $x = 1/1.84337061$, was found by power
iteration, and \eqref{eq:cw} and the strong connectivity were then verified in integer
arithmetic (\texttt{scripts/bounded\_certificate.py}).  Lemma~\ref{lem:cw} gives the claim.
\end{proof}

The theorem is computer-assisted: it relies on the program that constructs $\mathcal{A}_{12}$.
That program was checked against a brute-force count over all posets on at most $11$ points
(for $W \le 8$, in every cell of the table by number of elements and number of
join-irreducibles; see Section~\ref{sec:verify}).  Already $W = 6$, with $266$ states, gives
$\lambda_{\mathrm{vi}} \ge 1.7262$.

\boundedtable

Table~\ref{tab:bounded} shows the growth rates of $\mathcal{D}_W$.  They converge quickly, and
their limits agree with the ratios observed in Section~\ref{sec:dist}, which suggests
$\lambda_{\mathrm{vi}} = 1.72653630\ldots$ and $\beta = 1.84337061\ldots$\,.  Lower bounds
cannot show $\lambda_{\mathrm{vi}} < \beta$, the hypothesis of
Proposition~\ref{prop:growth}(b); for that an upper bound on $d_{\mathrm{vi}}$ below
$1.8433^n$ is needed, and the best known one is $2.33^n$~\cite{EHR}.

\section{Modular lattices}\label{sec:mod}

A finite lattice is modular if and only if it and its dual are semimodular.  Equivalently
it is graded with $\rk(x)+\rk(y) = \rk(x \wedge y) + \rk(x \vee y)$, and equivalently, for all
$x \ne y$,
\begin{equation}\label{eq:cover}
  x,y \text{ cover } x \wedge y \iff x \vee y \text{ covers } x \text{ and } y
\end{equation}
\cite[Section~3.3]{Stanley}.  There is no analogue of Birkhoff's theorem, and we work with
the lattice itself.

\subsection{Lower parts}
\begin{definition}
A \emph{modular lower part of height $r$} is a finite meet-semilattice $X$ with least
element, graded, all of whose maximal elements have rank $r$, such that for all
$x \ne y$ of equal rank $k$:
\begin{itemize}
\item[(U)] if $\rk(x \wedge y) = k-1$ and $k < r$, then $x$ and $y$ have a common upper cover;
\item[(D)] if $x$ and $y$ have a common upper cover, then $\rk(x \wedge y) = k-1$.
\end{itemize}
Its \emph{top layer} $T$ is the set of elements of rank $r$.  Elements $a \ne b$ of $T$ are
\emph{adjacent} if $\rk(a \wedge b) = r-1$.
\end{definition}

If $M$ is a finite modular lattice of rank at least $r$, then
$X_r(M) = \{x : \rk(x) \le r\}$ is a modular lower part of height $r$.  Conversely a modular lower part
with a single maximal element is a finite meet-semilattice with greatest element, hence a
lattice; a common upper cover of $x,y$ is then $x \vee y$, so (U) and (D)
give~\eqref{eq:cover}, and the lattice is modular.

Let $X$ be a modular lower part of height $r$.  Adjoin a set $T'$ of new elements, and
for each $c \in T'$ a nonempty set $D_c \subseteq T$ of elements declared to be covered
by $c$; let $X'$ be the resulting poset.

\begin{proposition}\label{prop:layer}
$X'$ is a modular lower part of height $r+1$ if and only if
\begin{itemize}
\item[(i)] every element of $T$ lies in some $D_c$;
\item[(ii)] distinct $a,b \in T$ are adjacent if and only if $\{a,b\} \subseteq D_c$ for some $c$;
\item[(iii)] for all $c \ne d$ in $T'$ the set
  $\{a \wedge b : a \in D_c,\ b \in D_d\}$ has a greatest element.
\end{itemize}
In that case the element in \textup{(iii)} is $c \wedge d$.
\end{proposition}

\begin{proof}
$X'$ is graded with all maximal elements of rank $r+1$ if and only if (i) holds.
Condition (ii) is (U) and (D) for pairs of rank $r$; for all other pairs (U) and (D) are
unchanged or vacuous.  It remains to show that $X'$ is a meet-semilattice if and only if
(iii) holds.

The common lower bounds of $c \ne d$ in $T'$ are the elements below $a \wedge b$ for
some $a \in D_c$, $b \in D_d$.  So $c \wedge d$ exists if and only if the set in (iii)
has a greatest element, which is then the meet.  Now assume (i) and (iii), and let
$c \in T'$, $x \in X$ with $x \not\le c$.  Choose $y \in T$ with $x \le y$; then
$y \notin D_c$, and by (i) there is $d \ne c$ with $y \in D_d$.  A common lower
bound $z$ of $c$ and $x$ satisfies $z \le c$ and $z \le d$, so $z \le c \wedge d =: e$,
and $e \in X$.  Hence the common lower bounds of $c$ and $x$ are those of $e$ and $x$,
and $c \wedge x = e \wedge x$ exists.
\end{proof}

Note that (iii) excludes two members $c \ne d$ with $|D_c \cap D_d| \ge 2$.  So in the
language of incidence geometry the sets $D_c$ of size at least two (``lines'') form a
partition of the edges of the adjacency graph of $T$ into cliques, and the elements with
$|D_c| = 1$ are new join-irreducibles.

\subsection{The frontier}
Let $Y(X) = \{a \wedge b : a,b \in T\} \supseteq T$ with the order induced from $X$, and let
$Y_1(X)$ be the set of its elements of rank $r-1$.  The \emph{frontier} of $X$ is the
triple $(Y,T,Y_1)$.

\begin{lemma}\label{lem:modfrontier}
Conditions \textup{(i)--(iii)} of Proposition~\ref{prop:layer} depend only on the
frontier of $X$ and the family $(D_c)_{c \in T'}$, and so does the frontier of $X'$.
\end{lemma}

\begin{proof}
Adjacency of $a,b$ means $a \wedge b \in Y_1$, and (iii) is a statement about the
order of $Y$.  For $X'$ we have $Y(X') = T' \cup \{c \wedge d : c \ne d\}$ where
$c \wedge d \in Y(X)$ is the element of (iii); an element $y \in Y(X)$ lies below $c$ if
and only if it lies below some $a \in D_c$; and $Y_1(X') = Y(X') \cap T$.
\end{proof}

\begin{lemma}\label{lem:modsym}
Let $\Gamma(X) \le \Sym(T)$ be the image of $\Aut(X)$.  Two families $(D_c)$,
$(D'_c)$ satisfying \textup{(i)--(iii)} give isomorphic extensions if and only if one
is mapped to the other, as a multiset of subsets of $T$, by an element of $\Gamma(X)$;
and $\Gamma(X')$ is determined by the frontier of $X$, $\Gamma(X)$ and the family.
\end{lemma}

The proof is the same as for Lemma~\ref{lem:sym}: isomorphisms preserve the rank, hence
restrict to lower parts.  As in Theorem~\ref{thm:dist} we obtain:

\begin{theorem}\label{thm:mod}
The number of isomorphism classes of modular lower parts of height $r+1$ with given
state $(Y,T,Y_1,\Gamma)$ and given number of elements satisfies the recursion of
Theorem~\ref{thm:dist}, with the $\Gamma$-orbits of families satisfying
\textup{(i)--(iii)} in place of birth vectors and $|T|$ as the weight of a state.
The modular lattices are the modular lower parts with $|T| = 1$, and they are vertically
indecomposable if and only if no layer strictly between bottom and top is a singleton.
\end{theorem}

\subsection{Representing the symmetry group}
Unlike distributive lattices, modular lattices contain large sets of interchangeable
elements at little cost: the lattice $M_k$ has $k+2$ elements and automorphism group
$\Sym(k)$.  We therefore store $\Gamma$ as a partition of $T$ into \emph{twin classes},
maximal sets on which $\Gamma$ contains the full symmetric group fixing all other points,
together with the explicit list $H$ of those elements of $\Gamma$ that are monotone on
every class.  Then $\Gamma = N H$ with $N$ the product of the symmetric groups of the
classes.  Membership of a given permutation of the new layer in $\Gamma(X')$ reduces, for each $h \in H$, to
comparing two multisets of incidence vectors, which gives the twin classes and the list
$H$ for $X'$.  Orbit representatives of families are found with a canonical form of the incidence
structure in which the twin classes are colour classes.

\subsection{Three bounds from modularity}
The local rules (i)--(iii) admit many lower parts that cannot be completed to a
modular lattice with at most $N$ elements.  The following consequences of modularity
remove most of them.  They are used only to discard states and families, so they affect
the running time but, if valid, not the result.  Let $M$ be a finite modular lattice and
$T$ its layer of rank $r$.

\begin{lemma}\label{lem:dist}
For $a,b \in T$ the number $\delta(a,b) = r - \rk(a \wedge b) = \rk(a \vee b) - r$ equals
the distance of $a$ and $b$ in the adjacency graph of $T$.
\end{lemma}

\begin{proof}
Let $m = a \wedge b$ and $\delta = \delta(a,b)$.  Take maximal chains
$m = w_0 \lessdot \dots \lessdot w_\delta = a$ and $m = u_0 \lessdot \dots \lessdot u_\delta = b$, and
put $c_i = w_{\delta-i} \vee u_i$.  Since $w_{\delta-i} \wedge u_i = m$, the rank formula gives
$\rk(c_i) = r$; also $c_0 = a$, $c_\delta = b$, and
$c_i \wedge c_{i+1} \ge w_{\delta-i-1} \vee u_i$, which has rank $r-1$.  So
consecutive $c_i$ are equal or adjacent, and the distance is at most $\delta$.
Conversely, if $a = c_0, c_1, \dots, c_p = b$ is a path, then
$\rk(c_0 \vee \dots \vee c_{j}) \le r + j$ by induction, because
$(c_0 \vee \dots \vee c_j) \wedge c_{j+1} \ge c_j \wedge c_{j+1}$ has rank at least $r-1$.
Hence $\delta \le p$.
\end{proof}

\begin{corollary}\label{cor:layers}
Let $\Delta$ be the diameter of the adjacency graph of $T$.  Then $M$ has rank at least
$r + \Delta$, and for $1 \le j < \Delta$ the layer of rank $r+j$ has at least as many
elements as any set of elements of $T$ with pairwise distances greater than $j$.  In
particular the adjacency graph is connected.
\end{corollary}

\begin{proof}
Elements $a,b$ with $\delta(a,b) > j$ have no common upper bound of rank $r+j$, and each of
them lies below some element of rank $r+j$.
\end{proof}

\begin{proposition}[Duality bound]\label{prop:dual}
For a multiset $\pi$ of sorted distance vectors, let $W(\pi)$ be the least number of
elements of rank less than $r'$ in any modular lattice $M'$ having a layer of rank $r'$
whose adjacency graph has distance profile $\pi$.  If the layer $T$ of $M$ has
profile $\pi$, then $M$ has at least $W(\pi)$ elements of rank greater than $r$.
\end{proposition}

\begin{proof}
The dual $M^*$ is modular, $T$ is one of its layers, and by Lemma~\ref{lem:dist} the
distances in the adjacency graphs of $T$ in $M$ and in $M^*$ are the same numbers
$\delta(a,b)$.  The elements of $M^*$ below $T$ are the elements of $M$ above $T$.
\end{proof}

The transfer matrix processes lower parts in increasing order of the number $a$ of
elements below the top layer.  When a state with profile $\pi$ is first created, at
$a = a_0$, we record $W_0(\pi) = a_0$.  A new state with $k$ elements in its top layer,
created at $a$, is discarded if $a + k + W_0(\pi) > N$, where $W_0(\pi)$ is read as
$a$ if $\pi$ has not been recorded.  By induction on $a$ this never discards the lower
part of a modular lattice $M$ with at most $N$ elements.  Indeed, let $c$ be the number of
elements of $M$ above the layer.  The lower part of $M^*$ with the same top layer has $c$ elements
below it and profile $\pi$, and it belongs to a modular lattice with at most $N$
elements.  If $c < a$, it has been created and kept, so $W_0(\pi) \le c$; if $c \ge a$,
then $W_0(\pi) \le a \le c$.  In both cases $a + k + W_0(\pi) \le |M| \le N$.  For
states with $a > N/2$ this means that a state survives only if its profile already
occurred in a sufficiently small lower part.  In our runs this rule and
Corollary~\ref{cor:layers} removed about three quarters of the states.

\begin{proposition}[Plane condition]\label{prop:plane}
Let $\ell_1 \ne \ell_2$ be elements of rank $r+1$ with $\ell_1 \wedge \ell_2 \in T$.  If
$m_1 \ne m_2$ of rank $r+1$ each cover two distinct elements of $T$ that lie below $\ell_1$ or
$\ell_2$, then $m_1 \wedge m_2 \in T$.
\end{proposition}

\begin{proof}
Both $m_i$ lie below $\ell_1 \vee \ell_2$, which has rank $r+2$.  So
$\rk(m_1 \wedge m_2) \ge 2(r+1) - (r+2) = r$.
\end{proof}

A violation of the plane condition would be detected by rule (iii) one layer later; testing
it while the family is built avoids enumerating families that cannot be extended.  The
clique partitions required by (ii) are found by backtracking over the uncovered adjacent
pairs.  Partial partitions are rejected up to the action of $\Gamma$ (equivalent partial
partitions have equivalent completions), and two counting bounds are applied to the
number of members still needed: the de~Bruijn--Erd\H{o}s theorem~\cite{dBE} (the members
meeting a set $Q$ of pairwise adjacent elements in at least two points form a linear space
on $Q$, so there are at least $|Q|$ of them unless one contains $Q$), and a pair-counting
bound.

\subsection{Results}
Table~\ref{tab:mod} contains the values for $n \ge 33$.  The run with $N=35$ took $19$
minutes on one core and $0.55$\,GB of memory, with $3.0\cdot 10^5$ states.  The runs with
$N=36$, $37$, $38$ took $37$, $73$ and $142$ minutes and used $4.5\cdot 10^5$,
$6.9\cdot 10^5$ and $1.06\cdot 10^6$ states; the cost grows by a factor of about $1.9$ per
element.  All values agree with those of~\cite{Koh21} for $n \le 35$ and
of~\cite[Table~1]{Koh23} for $36 \le n \le 38$, which were obtained by listing lattices or
racks, an entirely different method.  For comparison, the listing of the racks with $36$, $37$
and $38$ elements in~\cite{Koh23} took $1.6$, $3.6$ and $8.1$ core-days.  We have not
run $N = 39$ and $N = 40$.  The last
ratios are $m(\Nmod)/m(\Nmodprev) = \ratiom$ and
$m_{\mathrm{vi}}(\Nmod)/m_{\mathrm{vi}}(\Nmodprev) = \ratiomv$; they are still increasing,
as observed in~\cite{Koh21}.

\begin{table}[ht]
\tablemod
\caption{Values of $m(n)$ (A006981) and $m_{\mathrm{vi}}(n)$ (A342132) computed with the
transfer matrix.  They agree with~\cite{Koh21} ($n \le 35$) and~\cite{Koh23}
($n \le 38$).}\label{tab:mod}
\end{table}

\subsection{Growth rate}
The construction of Section~\ref{sec:bounded} applies to modular lattices as well, and here
finiteness is clear: if every layer has at most $W$ elements, a frontier $Y$ has at most
$W + \binom{W}{2}$ elements, so there are finitely many states.  Let $\mathcal{A}^{\mathrm{mod}}_W$
be the automaton of the modular lattices all of whose layers have at most $W$ elements
(vertices: the states reachable from the one-element layer; weights: $|T|$).

\begin{theorem}\label{thm:modgrowth}
$\limsup_n m_{\mathrm{vi}}(n)^{1/n} \ge 2.3951$ and $\lim_n m(n)^{1/n} \ge 2.4182$.  More precisely
$m_{\mathrm{vi}}(n) = \Omega(c^n)$ for every $c < 2.3951$.
\end{theorem}

\begin{proof}
As for Theorem~\ref{thm:growth}, with $\mathcal{A}^{\mathrm{mod}}_8$, which has $31\,435$
vertices, and its strong components with $29\,740$ (without the one-element layer) and
$29\,742$ vertices.  The limit for $m$ exists in $(0,\infty]$ because $m(n+1)$ is
supermultiplicative.  The components have period one, i.e.\ the weights of their closed
walks have greatest common divisor one; the numbers $c_s(n)$ of Lemma~\ref{lem:cw} are
supermultiplicative and positive for all large $n$, so $c_s(n)^{1/n}$ converges to its
supremum, which gives the last statement.
\end{proof}

The previous lower bound was $\Omega(2.3122^n)$~\cite{Koh21}; no exponential upper bound is
known.  Table~\ref{tab:boundedmod} shows the growth rates for $W \le 8$.  Unlike in the
distributive case they have not converged; their increments ($0.042$, $0.027$ for
$W = 7, 8$) decrease, but we cannot tell from these data whether they stay bounded.  As a
check of the two automaton programs, $\mathcal{A}^{\mathrm{mod}}_W$ restricted to distributive
lattices has the same growth rates as $\mathcal{A}_W$ for $W \le 6$, to all ten computed digits,
although its states are different.

\boundedtablemod

\section{Semimodular and geometric lattices}\label{sec:other}

The proof of Proposition~\ref{prop:layer} treats its three conditions separately:
\textup{(i)} makes $X'$ graded, \textup{(iii)} makes it a meet-semilattice, and the two
directions of \textup{(ii)} are the rules (U) and (D) for the pairs of rank $r$.  So the
frontier $(Y,T,Y_1)$ of Section~\ref{sec:mod} serves, without change, for the lower parts of
\begin{itemize}
\item graded lattices: conditions (i) and (iii);
\item (upper) semimodular lattices: (i), (iii) and ``adjacent elements lie in a common $D_c$'';
\item lower semimodular lattices: (i), (iii) and ``the elements of each $D_c$ are pairwise adjacent'';
\item geometric lattices, i.e.\ semimodular lattices in which every element is a join of atoms:
      as for semimodular lattices, with $|D_c| \ge 2$ for all $c$ of rank at least two.
\end{itemize}
Lemmas~\ref{lem:modfrontier} and~\ref{lem:modsym} and Theorem~\ref{thm:mod} hold for each of
these classes with the same proofs.  By (iii) two members of a family still share at most
one element, so the members with at least two elements form a partial linear space on $T$
which, in the semimodular case, covers all adjacent pairs and may cover others.

The bounds of Section~\ref{sec:mod} that rest on modularity are not available.  In fact a
semimodular lower part can always be completed by a single top element, so no state can be
discarded, and the lattices themselves are far more numerous than modular ones.
Nevertheless the number of states stays moderate when the lattices are built from the bottom:
for $N = 24$ there are $8.9 \cdot 10^{5}$ states for about $1.5 \cdot 10^{9}$ semimodular
lattices.  The direction matters.  Counting the dual class (lower semimodular lattices from the
bottom) gives the same numbers, which we used as a check for $n \le 16$, but needs vastly more
states.  For graded lattices in general the enumeration of the families is too expensive for
the method to compete with generation~\cite{Koh19}.

The run with $N = 26$ took $69$ minutes on one core and $7.9$\,GB of memory, with
$5.9\cdot 10^{6}$ states.  For $N = 27$ we used a version of the program that stores each
state only as its canonical key, does not store transitions (a state is expanded again
whenever it is reached with a new weight, about $1.6$ times on average), and expands the
states of one weight level on four threads; it took $95$ minutes and $9.2$\,GB, with
$1.54\cdot 10^{7}$ states.  The runs reproduce the values for $n \le 25$ of~\cite{Koh19} and
give two new values (Table~\ref{tab:other}).  Memory is the limit here.  The value for $n=27$
is computed at $n = N$; since no state is ever discarded in this class, the bound $N$ only
truncates the enumeration of families.

For geometric lattices (simple matroids counted by their number of flats) the same program
reproduces the known values for $n \le 35$, but it offers no advantage: these lattices have
small rank, the frontier is essentially the whole lattice, and the number of states is close
to the number of lattices.  The cost is then dominated by the enumeration of linear spaces
on the atoms, and we did not extend that sequence.

\resultsother

\section{Sequences determined by the tables}\label{sec:further}

Two further OEIS sequences are functions of $d$ and $d_{\mathrm{vi}}$.  The numbers $w_k$ of
canonically $2$-indecomposable posets with $k$ antichains (A072407) satisfy
$d_{\mathrm{vi}}(k) = \sum_{j=0}^{k-1} w_{k-j}\, d_{\mathrm{vi}}(j)$ for $k > 1$, with the conventions
$d_{\mathrm{vi}}(0) = 1$, $d_{\mathrm{vi}}(1) = w_1 = 0$~\cite[Corollary~7]{EHR}; this gives $w_k$ for
$k \le \Ndist$ (previously $k \le 49$).  The numbers $n \le \Ndist$ for which $d(n)$ is prime
(A186808) are $4, 5, 6, 10, 12, 13, 18, 21, 23, 26, 66, 88$; the last two are new, and
primality was proved with PARI/GP.

\section{Verification}\label{sec:verify}

The results are the output of a few C++ programs of a few hundred lines each.  We
list what supports them, and what does not.

\begin{enumerate}
\item \emph{Known values.}  The distributive program reproduces $d(n)$ for $n \le 60$
and $d_{\mathrm{vi}}(n)$ for $n \le 49$ as given in the OEIS, and $d_{\mathrm{vi}}(n)$
for $50 \le n \le 60$ as printed in~\cite[Table~2]{Koh21}.  The modular
program reproduces $m(n)$ and $m_{\mathrm{vi}}(n)$ for $n \le 35$~\cite{Koh21} and for
$36 \le n \le 38$~\cite{Koh23}.  The earlier
computations are structurally different from ours.
\item \emph{Brute force beyond the known range.}  We generated all $49\,519\,391$ posets on at
most $11$ points with \texttt{genposetg} from the \texttt{nauty} package and tabulated those whose down-set lattice
is vertically indecomposable by (number of down-sets, number of points).  The
distributive transfer matrix, refined by the number of layers, agrees in all $443$ cells
with at most $110$ down-sets; $192$ of these cells, with $25\,179\,538$ lattices, have
more than $60$ down-sets.
\item \emph{Planar distributive lattices.}  Restricted to posets of width at most two the
distributive automaton has $23$ states for $N = 160$.  Its output agrees for all
$n \le 160$ with a direct count of the regions between two lattice paths, and with
OEIS sequence A343161.
\item \emph{Duality.}  Both transfer matrices run from the bottom and are not symmetric
under duality.  The number of vertically indecomposable lattices with $n$ elements, $s$
atoms and $t$ coatoms must be symmetric in $(s,t)$.  This holds in all $3013$ cells for
distributive lattices with $n \le 160$ and in all cells for modular lattices up to the
largest $N$.
\item \emph{Two programs on distributive lattices.}  A modular lattice of finite length
is distributive if and only if no interval of length two has three atoms.  With this
restriction the modular program gives $d_{\mathrm{vi}}(n)$, and it agrees with the
distributive program for all $n \le 40$.  The two programs share no state definition,
symmetry handling or pruning rule.
\item \emph{Different bounds $N$.}  All pruning depends on $N$, and the value for
$n = N$ is computed with no slack.  Runs with different $N$ agree on their common range
(distributive: $N = 120, 130, 150, 160$, the first two also with an independent reference
implementation; modular: $N = 26, 28, 30, 32, 35, 36, 37, 38$; semimodular: $N = 24, 25, 26, 27$).  The largest
value in the tables for distributive and for semimodular lattices has no such confirmation from
a run with larger $N$.  For modular lattices every boundary value $m(N)$, $N \le 38$, was
correct.
\end{enumerate}

For the automata of Section~\ref{sec:bounded} we compared, for $W = 3,\dots,8$, the numbers of
walks with at most $11$ steps and weight at most $119$ with a brute-force count of the
vertically indecomposable lattices $J(P)$, $|P| \le 11$, with all layers of size at most $W$;
they agree in every cell (up to $10\,340\,439$ lattices for $W = 8$).

For semimodular lattices the checks are the known values for $n \le 25$, the runs with
different $N$, and the agreement with the count of the dual class for small $n$.

There is no independent method for distributive lattices of width at least three with more than
$110$ elements, nor for semimodular lattices with more than $25$ elements; for
these the evidence is items (4)--(6).  An independent recomputation would be valuable.

\section{Remarks}

\begin{remark}
We first tried to count posets by their number of down-sets through modular
decomposition: the number of down-sets of a lexicographic sum over a poset $Q$ is
$\sum_{A} \prod_{x \in A} (i_x - 1)$, the sum being over the antichains $A$ of $Q$ and $i_x$
denoting the number of down-sets of the $x$-th summand.  So it suffices to generate the prime
posets.  But the number of prime posets with $n$ down-sets still appears to grow like
$1.64^n$ (for $n \le 36$), so nothing is gained over~\cite{Koh21}.
\end{remark}

\begin{remark}
The method needs a grading that is respected by isomorphisms and a bounded description
of what the upper part can see of the lower part.  It does not apply as it stands to
all lattices, which are not graded.
\end{remark}

\begin{remark}
Memory, not time, limits the distributive computation: at $N = 160$ there are
$1.75 \cdot 10^7$ states.  For modular lattices the enumeration of clique partitions for
frontiers with many pairwise adjacent elements dominates the running time.
\end{remark}

\subsection*{Declaration on the use of AI}
This work was carried out with substantial assistance from an AI system (Claude, by
Anthropic), used interactively under the direction of the author.  The AI system proposed and
implemented the algorithms, wrote the programs, ran the computations, and drafted the proofs
and the text.  The author has reviewed the paper and takes full responsibility for its
content.  The numerical results rest on the checks described in Section~\ref{sec:verify}.

\subsection*{Data availability}
The programs, the logs of all runs and the complete tables are available at~\cite{code}
and as ancillary files.


\begin{thebibliography}{99}

\bibitem{BM} G.~Brinkmann and B.~D.~McKay, Posets on up to 16 points,
  \emph{Order} \textbf{19} (2002), 147--179.

\bibitem{dBE} N.~G.~de~Bruijn and P.~Erd\H{o}s, On a combinatorial problem,
  \emph{Indag.\ Math.} \textbf{10} (1948), 421--423.

\bibitem{EHR} M.~Ern\'e, J.~Heitzig and J.~Reinhold, On the number of distributive
  lattices, \emph{Electron.\ J.\ Combin.} \textbf{9} (2002), \#R24.

\bibitem{FS} P.~Flajolet and R.~Sedgewick, \emph{Analytic Combinatorics}, Cambridge
  University Press, 2009.

\bibitem{GT} V.~Gebhardt and S.~Tawn, Constructing unlabelled lattices,
  \emph{J.\ Algebra} \textbf{545} (2020), 213--236.

\bibitem{HR} J.~Heitzig and J.~Reinhold, Counting finite lattices,
  \emph{Algebra Universalis} \textbf{48} (2002), 43--53.

\bibitem{JL} P.~Jipsen and N.~Lawless, Generating all finite modular lattices of a
  given size, \emph{Algebra Universalis} \textbf{74} (2015), 253--264.

\bibitem{Koh19} J.~Kohonen, Generating modular lattices of up to 30 elements,
  \emph{Order} \textbf{36} (2019), 423--435.

\bibitem{Koh21} J.~Kohonen, Cartesian lattice counting by the vertical 2-sum,
  \emph{Order} \textbf{39} (2022), 113--141.

\bibitem{Koh23} J.~Kohonen, Simplifying modular lattices by removing doubly irreducible
  elements, preprint, arXiv:2311.16643, 2023.

\bibitem{nauty} B.~D.~McKay and A.~Piperno, Practical graph isomorphism, II,
  \emph{J.\ Symbolic Comput.} \textbf{60} (2014), 94--112.

\bibitem{oeis} OEIS Foundation Inc., \emph{The On-Line Encyclopedia of Integer
  Sequences}, \url{https://oeis.org}.

\bibitem{Stanley} R.~P.~Stanley, \emph{Enumerative Combinatorics}, Vol.~1, 2nd ed.,
  Cambridge University Press, 2012.

\bibitem{code} Programs and data for this paper,
  \url{https://github.com/mensch72/counting}.

\end{thebibliography}
\end{document}